\documentclass[11pt]{article}
\usepackage{graphicx}
\newtheorem{theorem}{Theorem}
\newtheorem{lemma}[theorem]{Lemma}

\newtheorem{proposition}[theorem]{Proposition}

\newtheorem{corollary}[theorem]{Corollary}

\newtheorem{definition}[theorem]{Definition}

\newtheorem{notation}[theorem]{Notation}

\begin{document}

\title{Resolutions of Segre embeddings of projective spaces of any dimension}
\author{  Elena Rubei }
\date{\hspace*{2cm}}
\maketitle

\vspace{-1cm}

{\small {\bf Abstract.} 
This paper deals with syzygies  of the ideals of the Segre embeddings.
Let $d \geq 3 $ and $n_1,..., n_d \in {\bf N} - \{0\}$.
 We prove  that  ${\cal O}_{{\bf P}^{n_1} \times .... \times 
{\bf P}^{n_d}}(1,...,1)$   satisfies Green-Lazarsfeld's Property $N_p$ 
 if and only if $p  \leq  3$.   }

\bigskip

\section{Introduction}

Let  $ L $ be a very ample  line bundle on a smooth
complex projective variety $Y$ and let  $\varphi_{L}: Y \rightarrow
{\bf P}(H^{0}(Y, L)^{\ast})$ be  the map  associated to $L$.
We  recall the definition of Property $N_{p}$ of Green-Lazarsfeld,
 studied for the first time  by Green in  \cite{Green1}
(see also \cite{G-L}, \cite{Green2}):

{\em
let $Y$ be a smooth complex projective variety and let $L$ be a very ample
 line bundle on $Y$ defining an embedding $\varphi_{L}: Y \hookrightarrow
{\bf P}={\bf P}(H^{0}(Y,L)^{\ast })$;
set $S= S(L)= \oplus_n Sym^{n} H^{0}(L) $,
the homogeneous coordinate ring of the projective space ${\bf P}$, and
consider
the graded $S$-module  $G=G(L)= \oplus_{n} H^{0}(Y, L^{n})$; let $E_{\ast}$
\[ 0 \longrightarrow E_{l} \longrightarrow E_{l-1} \longrightarrow 
... \longrightarrow E_{0} \longrightarrow G \longrightarrow 0\] 
be a minimal graded free resolution of $G$;
the line bundle $L$ satisfies Property
$ N_{p}$ ($p \in {\bf N}$) iff

\hspace{1cm} $E_{0}= S$

\hspace{1cm} $E_{i}= \oplus S(-i-1)$    \hspace{1cm}
for $1 \leq i \leq p $.}

(Thus $L$ satisfies Property $N_{0}$  iff  $Y \subset
{\bf P}(H^{0}(L)^{\ast })$ is projectively normal;
 $L$ satisfies Property $N_{1}$
iff $L$ satisfies $N_{0}$
and the homogeneous ideal $I$ of $Y \subset
{\bf P}(H^{0}(L)^{\ast })$ is generated by quadrics;
$L$ satisfies  $N_{2}$  iff
 $L$ satisfies  $N_{1}$ and the module of syzygies among
quadratic
generators $Q_{i} \in I$ is spanned by relations of the form
 $\sum L_{i}Q_{i}=0$, where $L_{i}$  are linear polynomials;
and so on.)

In this paper we will consider the case of Segre embedding i.e. 
 the case  $Y= {\bf P}^{n_1} \times ..... \times {\bf P}^{n_d} $,
 $L ={\cal O}(1,....,1)$.

Among the papers on syzygies of Segre-Verones embeddings
 we quote  \cite{B-M}, \cite{Green1}, \cite{O-P},
\cite{J-P-W}, \cite{G-P}, \cite{Las}, \cite{P-W}, \cite{Ru3}, \cite{Ru4}, 
\cite{S-S}. 

In the first one, the authors examined the cases in which the resolution 
is ``pure'' i.e. the minimal generators of each module of syzygies have 
the same degree. We quote the following 
 results from the other papers:

{\it Case $d=1$}, i.e. the case of the Veronese embedding:

\begin{theorem} \label{Gr} ({\bf Green}) 
\cite{Green1}.
Let $a$ be a positive integer.
The line bundle ${\cal O}_{{\bf P}^{n}} (a)$ 
  satisfies Property $N_{a +1}$.
\end{theorem}

\begin{theorem} \label{SS}  ({\bf Schenck-Smith}) 
 \cite{S-S}.
Let $a$ be a positive integer.
The line bundle ${\cal O}_{{\bf P}^{n}} (a)$ 
  satisfies Property $N_{a +1}$.
\end{theorem}

 The case ${\cal O} (3) $ satisfies $N_4$ is also proved by {\bf Rubei} in 
\cite{Ru4}.

\begin{theorem} \label{OP}  ({\bf Ottaviani-Paoletti}) \cite{O-P}.
If $n \geq 2$,
  $a \geq 3$ and the bundle
${\cal O}_{ {\bf P}^{n} } (a)$ satisfies Property $N_{p} $, then $
p \leq 3a -3$.
\end{theorem}

\begin{theorem} \label{JPW} ({\bf Josefiak-Pragacz-Weyman}) \cite{J-P-W}.
The bundle
${\cal O}_{ {\bf P}^{n} } (2)$ satisfies Property $N_{p} $ if and only if $
p \leq 5$ when $n \geq 3$ and  for all $p$ when $n=2$.
\end{theorem}

(See \cite{O-P} for a more complete bibliography.)

{\it
Case $d=2$:}

\begin{theorem} \label{GP} ({\bf Gallego-Purnapranja}) \cite{G-P}.
Let $a,b \geq 2$. 
The line bundle ${\cal O}_{{\bf P}^{1} \times {\bf P}^{1} }(a,b)$ 
 satisfies Property $N_p$ if and only if $p \leq 2a+2b-3 $.
\end{theorem}

\begin{theorem} \label{LPW} ({\bf Lascoux-Pragacz-Weymann}) \cite{Las}, 
\cite{P-W}. Let $n_1, n_2 \geq 2$.
   The line bundle ${\cal O}_{{\bf P}^{n_1} \times {\bf P}^{n_2}}(1,1)$ 
 satisfies Property $N_p$ if and only if $p \leq 3$.
\end{theorem}

As to the case $d \geq 3$, in \cite{Ru3} we proved

\begin{theorem} ({\bf Rubei}) \label{P1}.
Let $ d \geq 3 $. The line bundle ${\cal O}_{{\bf P}^{1} 
\times ...\times {\bf P}^{1} }(1,..., 1)$ ($d$ times)
 satisfies Property $N_p$ if and only if $p \leq 3$.
\end{theorem}

Here we prove 

\begin{theorem} \label{pasqua}.
Let $ d \geq 3 $ and $n_1,..., n_d \in {\bf N}-\{0\}$.
 The line bundle ${\cal O}_{{\bf P}^{n_1} 
\times ...\times {\bf P}^{n_d} }(1,..., 1)$ 
 satisfies Property $N_p$ if and only if $p \leq 3$.
\end{theorem}

To prove it we use  Theorem \ref{P1}. As in the proof of it, we use a 
topological approach, but the problems that arise
in proving  Theorem \ref{pasqua} 
  are more similiar to the 
problems that arise if we try to solve the open problem of the syzygies of 
the Veronese embeddings.

\section{A preliminary proposition}

\begin{proposition} \label{pasquetta}.
Let $ d \geq 3 $ and $n_1,..., n_d \in {\bf N}$ with $n_i \geq p$ for some $i$.
 The line bundle ${\cal O}_{{\bf P}^{n_1} 
\times ...\times {\bf P}^{n_d} }(1,..., 1)$ 
 satisfies $N_p$ if and only if ${\cal O}_{{\bf P}^{n_1} 
\times ...\times {\bf P}^{n_{i-1}}\times {\bf P}^{p} \times 
 {\bf P}^{n_{i+1}} \times ...\times  {\bf P}^{n_{d}}}(1,..., 1)$ 
 satisfies  $N_p$.
\end{proposition}

{\em Proof.}
 Let   $ L $ be a very ample  line bundle on a smooth
complex projective variety $Y$. We recall from \cite{Green2} that $L$
satisfies $N_p$ iff $$ (Tor^{S(L)}_p (G(L),{\bf C}) 
)_{p+q}= 0 \;\; \forall q \geq 2$$
(see Introduction for the notation) and
$ (Tor^{S(L)}_p (G(L),{\bf C}) )_{p+q}$ is equal to the homology of
the Koszul complex
$$\wedge^{p+1} H^0 (L)  \otimes H^0(L^{q-1})  \rightarrow 
\wedge^{p} H^0 (L)  \otimes H^0(L^{q})
\rightarrow \wedge^{p-1} H^0 (L)  \otimes H^0(L^{q+1})$$

Now let $Y=  {\bf P}(V_1) \times .... \times {\bf P} (V_d)$
 and $L={\cal O}(1,...., 1)$. In our case   
$ (Tor^{S(L)}_p (G(L),{\bf C}))_{p+q}$ is equal to the homology of
$$ (\ast_{V_1 , ...,V_d}) \;\;\;\;\;\;\;\;\;\;\;\; \begin{array}{l}
 \wedge^{p+1}  (V_1 \otimes .... \otimes V_d) \otimes Sym^{(q-1)} V_1 
 \otimes ..... \otimes Sym^{(q-1)} V_d ;\;\;\;\;\;\;
\stackrel{\alpha^{V_1 \times .... \times V_d}_{p+1,q-1}}{\longrightarrow} \\ 
\wedge^{p} (V_1 \otimes .... \otimes V_d) \otimes Sym^{q} V_1 
 \otimes ..... \otimes Sym^{q} V_d ;\;\;\;\;\;\;
\stackrel{\alpha^{V_1 \times .... \times V_d}_{p,q}}{\longrightarrow}  \\
\wedge^{p-1} (V_1 \otimes .... \otimes V_d) \otimes Sym^{(q+1)} V_1 
 \otimes ..... \otimes Sym^{(q+1)} V_d ;\;\;\;\;\;\;
\stackrel{\alpha^{V_1 \times .... \times V_d}_{p-1,q+1}}{\longrightarrow} 
\end{array}
$$ 
since the maps are $GL(V_i)$-invariant, 
 $ (Tor^{S(L)}_p (G(L),{\bf C}) )_{p+q}$
is a $GL(V_i)$-module,
 see Rem. of \S 2 of \cite{Green1} and Prop. 1.8
\cite{O-P}. 
The Young diagrams of the irreducible subrepresentations of
the $GL(V_i)$-module 
$\wedge^{p} (V_1 \otimes .... \otimes V_d) \otimes Sym^{q} V_1 
 \otimes ..... \otimes Sym^{q} V_d  $ have  at most $p+1 $ rows (see for 
instance p.80 \cite{F-H}).  Thus the Young diagrams of the
 irreducible subrepresentation of $ (Tor^{S(L)}_p (G(L),{\bf C}) )_{p+q}$ 
have at most  $p  +1 $ rows and one can easily prove 
that  these Young diagrams are the same for  $V_1,...,V_d$ and for 
$V_1,....,V_{i-1},V'_i, V_{i+1},...., V_d$ with 
$dim V'_i =p+1$ 
(by considering  an injective map $V'_i \rightarrow V_i$ and the  induced 
commutative diagram with second row  $(\ast_{V_1 , ...,V_d})$ and first 
row  $(\ast_{V_1 , ....,V_{i-1},V'_i,V_{i+1}, ...,V_d})$ seeing every term
of the two rows  as $Gl(V_i)$ and $GL(V'_i)$ representations respectively).

 Thus,
since  the Young diagrams of the
 irreducible subrepresentations of $ (Tor^{S(L)}_p (G(L),{\bf C}) )_{p+q}$ 
have at most  $p  +1 $ rows,
we have that if these representations  are zero for $dim(V_i) = 
p + 1 $  they are zero also for $dim(V_i) \geq p  +1 $.

\hfill \framebox(7,7)

\section{Recalls on syzygies of toric ideals and on syzygies on the product of
 two varieties}

We  recall   some facts on toric ideals from \cite{St}.
Let $k \in {\bf N}$. Let $A = \{a_1,..., a_m \} \subset {\bf Z}^k$.
The toric ideal  ${\cal I}_A$ is defined as  the  ideal 
in ${\bf C}[x_1, ..., x_m]$
generated as vector space by the
binomials 
$$x_1^{u_1}... x_m^{u_m}
- x_1^{v_1}... x_m^{v_m}$$
 for $(u_1,..., u_m),(v_1,.., v_m) \in {\bf N}^m$, 
with $\sum_{i=1,...,m} u_i a_i = \sum_{i=1,...,m} v_i a_i $. 
           
We have that ${\cal I}_A$ is homogeneous iff $\exists \;
\omega \in {\bf Q}^k$ s.t. $\omega  \cdot a_i= 1$ $\forall i = 1,..., m$;
the rings ${\bf C}[x_1, ..., x_m]$ and 
  ${\bf C}[x_1, ..., x_m]/{\cal I}_A$ 
are multigraded by ${\bf N}A$ via $\deg x_i =a_i$;  
the element $x_1^{u_1}...x_m^{u_m}$ 
has multidegree $ b =  \sum_i u_i a_i \in {\bf N}A$ and degree 
$\sum_i  u_i = b \cdot \omega$; we define $\deg b = b \cdot \omega$.

For each $b \in {\bf N}A$, let 
$\Delta_b$ be the simplicial complex (see \cite{Sp})
on the set $A$  defined   as follows:
 $$ \Delta_b = \{ \langle F \rangle 
 | \; F \subset A : b - \sum_{a \in F} 
a  \in {\bf N}A \}$$

The following theorem studies the syzygies of the
ideal ${\cal I}_A$; it was proved by Campillo and Marijuan for $k=1$ 
in \cite{C-M} and by Campillo and Pison for general $k$ and $j=0$ in 
\cite{C-P}; the following more general statement is due to Sturmfels  
(Theorem 12.12 p.120 in \cite{St}).

\begin{theorem} \label{CPS}
({\rm see \cite{St}, \cite{C-M}, \cite{C-P}}).
Let $A=\{a_1,...,a_m\} \subset {\bf N}^k$
 and ${\cal I}_A$ be the associated 
toric ideal.  Let $ 0 \rightarrow E_n \rightarrow ...\rightarrow
 E_1 \rightarrow E_0 \rightarrow G \rightarrow 0$ be a 
minimal free resolution of $G= {\bf C}[x_1,..., x_m]/ {\cal I}_A$
on ${\bf C}[x_1,..., x_m]$.
Each of the generators of $E_j$ has a unique multidegree. 
The number of the generators of 
multidegree $b \in {\bf N}A$ of  $E_{j+1}$  equals the rank of the 
$j$-th
reduced homology group $\tilde{H}_j(\Delta_b, {\bf C}) $. 
\end{theorem}

Finally we recall the following proposition 
proved in \cite{Ru3}

\begin{proposition} ({\bf Rubei}) \label{prodotti}.
Let $X$ and $Y$ be two projective varieties and let
$L$ be a line bundle on $X$ and $M$ a line bundle on $Y$.
Let $\pi_X : X \times Y \rightarrow X$ and $\pi_Y : X \times Y \rightarrow Y$
be the canonical projections. 
Suppose $L$ and $M$ satisfy Property $N_1$. Let $p \geq 2$.
If $L$ does not satisfy Property $N_p$, then  $\pi_{X}^{\ast}L
\otimes \pi_{Y}^{\ast}M$ does not satisfy Property $N_p$, either. 
\end{proposition}

\section{Proof of Theorem \ref{pasqua}}

\begin{notation}
$\bullet$ homologous means homologous in the reduced homology.

$\bullet$ $\sim_A$ means homologous in $A$, i.e. $\gamma \sim_A \gamma'$
means that $\exists \, \beta$ chain in $A$ s.t. $\partial \beta = 
\gamma - \gamma'$.

$\bullet$ $e_i$ denotes the $i$-th 
element of the canonical basis of ${\bf R}^n$.

$\bullet $ The symbol $\ast $ denotes the joining.

$\bullet $ For any  $v \in {\bf R}^n$ 
$v_i$ denotes the $i$-th coordinate, that is the lower index denotes the 
coordinate. 
\end{notation}

If we take $A =A_{n_1,..., n_d}= \{ 
(x^0_1, ...,
x^{n_1}_{1},..........., x^0_d, ...,
x^{n_d}_{d})
| \sum_{j=1 ,..., n_i}  x^j_i = 1\;\;  x^j_i \in {\bf N}\}$, 
     we have that ${\cal I}_{A_{n_1,..., n_d}}$ is the ideal of the  
embedding of ${\bf P}^{n_1} \times ...\times {\bf P}^{n_d} $
 by  ${\cal O}(1,...,1)$.
In this case  $\omega =\omega_d = \frac{1}{d} (1, ......, 1)$.

Let $b \in {\bf N} A_{n_1,..., n_d}$; 
we have that  $\deg b=(= b \cdot \omega) =k$ iff
$b$ is the sum of $k$ (not necessarily distinct) elements of 
$A_{n_1,..., n_d}$.
Observe that a simplex $S$ 
 with vertices in $A_{n_1,..., n_d}$ is a simplex of $\Delta_b$
iff  the sum $s$ of the vertices  of $S$ is s.t.  $s_i \leq b_i$ 
$\forall i =1,....,n+1$.
We generalize the definition of the simplicial complex 
$\Delta_b$ in \S 3   in the following way:

\begin{notation} Let $v \in {\bf N}^{n+1}$. Let $\Delta_v$ be  
the following simplicial complex: a
simplex $S$ with vertices in $A_{n_1,..., n_d}$ is a simplex of $\Delta_v$
iff  the sum $s$ of the vertices  of $S$ is s.t.  $s_i \leq v_i$ 
$\forall i =1,....,n+1$.
\end{notation}

The main points of the  proof of the Thm. \ref{pasqua}
  are Propositions \ref{step1} and \ref{step2}.

\begin{notation} Let $d,n_1,...,n_d  \in {\bf N}-\{0\}$ and 
$b \in {\bf N} A_{n_1,..., n_d}$. Let $X_b$ 
be the following simplicial complex on $A_{n_1,..., n_d}$: 
$$
X_b := \Delta_b \cup \Delta_{
\tiny{
\left( \hspace{-0.2cm} \begin{array}{c}
b_1 -1\\
b_2 +1 \\
b_3 \\
. \\
.\\
.
\end{array}
 \hspace{-0.2cm}
\right)}} \cup .... \cup \Delta_{
\tiny{
\left( \hspace{-0.2cm} \begin{array}{c}
0\\
b_2 + b_1 \\
b_3 \\
. \\
.\\
.
\end{array}
 \hspace{-0.2cm}
\right)}} $$ (in the obvious sense that a simplex with vertices in
 $A_{n_1,..., n_d}$  is a simplex of $X_b$ iff  it is a simplex of 
$\Delta_{b - k e_1 +k e_2}$ for some $k \in \{0,..., b_1\}$)
\end{notation}

\begin{proposition} \label{step1} Let $d,n_1,...n_d, p \in {\bf N} -\{0\}$.
Let $b \in {\bf N} A_{n_1,..., n_d}$ with $deg(b) \geq p+2$.
Let $\gamma $ be a $(p-1)$-cycle in $\Delta_b$.
If the following conditions hold

a) $\tilde{H}_{p-3}(\Delta_{c- e_1  }) =0 $
 $\forall c  \in {\bf N} A_{n_1,..., n_d}$ with 
  $deg(c) \geq p+1$  if $p \geq 3$

b) ${\cal O}_{{\bf P}^{n_1} \times ....\times {\bf P}^{n_d} 
}(1,...,1)$ satisfies Property $N_{p-1}$,

then  $\exists \gamma'$ cycle  in $\Delta_{
\tiny{
\left( \hspace{-0.2cm} \begin{array}{c}
0\\
b_1+b_2 \\
b_3 \\
. \\
.\\
.
\end{array}
 \hspace{-0.2cm}
\right)}}$ s.t. $\gamma \sim_{X_b} \gamma'$. 
\end{proposition}

\begin{proposition} \label{step2} Let $d,n_1,...n_d \in {\bf N} -\{0\}$ and
 $p \in \{  2,3\}$.
Let $b \in {\bf N} A_{n_1,...,n_d}$  with $deg(b) \geq p+2$.
Let $\gamma $ be a $(p-1)$-cycle in $\Delta_b$. If   $\gamma \sim_{X_b} 0$
 then $ \gamma \sim_{\Delta_b} 0 $.
\end{proposition}

\begin{lemma} \label{HOc-ei}  Let $d,n_1,...,n_d  \in {\bf N} -\{0\}$ and
 $c \in {\bf N} A_{n_1,..., n_d}$ with $deg (c) \geq 4$.
 Then $\tilde{H}_0(\Delta_{c-e_i})=0 $ $\forall i$.
\end{lemma}

\smallskip

We show now how  Thm. \ref{pasqua} follows  from Propositions
 \ref{step1} and  \ref{step2} and  Lemma \ref{HOc-ei}.

\bigskip

{\it Proof of Thm. \ref{pasqua}}. By Prop. \ref{pasquetta}, 
 to prove ${\cal O}_{{\bf P}^{n_1} \times ...\times {\bf P}^{n_d}}(1,...,1)$
  satisfies $N_3$, 
 it is sufficient to prove our statement when $n_i  \leq 3$ $\forall i$.
 We prove the 
statement by induction on $\sum_i n_i$; the first step of induction is given
by Thm. \ref{P1}.   
By Thm. \ref{CPS}, the bundle 
 ${\cal O}_{{\bf P}^{n_1} \times ...\times {\bf P}^{n_d}}(1,...,1)$
  satisfies $N_p$ iff 
 $\tilde{H}_{q-1}(\Delta_b)=0$  $\forall b \in {\bf N}A_{n_1,..., n_d}$ with 
$\deg b  \geq q+2$ $\forall q \leq p$;
in particular,  in order to prove that  it
 satisfies $N_3$, we have to prove 
that $\tilde{H}_{p-1}(\Delta_b)=0$  
$ \forall b \in {\bf N}A_{n_1,..., n_d}$ with 
$\deg b  \geq p+2$ for $p=2,3$.

We show that {\em ($\star$) if ${\cal O}_{{\bf P}^{n_1} \times ...\times 
{\bf P}^{n_d}}(1,...,1)$  satisfies $N_{p-1}$ 
then $\tilde{H}_{p-1}(\Delta_b)=0$  
$ \forall b \in {\bf N}A_{n_1,..., n_d}$ with 
$\deg b  \geq p+2$ for $p \in\{ 2,3\} $. }  
Let $b \in {\bf N} A_{n_1,..., n_d}$ with $deg(b) \geq p+2$, $p \in \{2,3 \}$.
Let $\gamma $ be a $(p-1)$-cycle in $\Delta_b$. We want to prove  
$\gamma \sim_{\Delta_b} 0$.  
 By   Prop. \ref{step1}, $ \exists \gamma'$ 
cycle   in $\Delta_{
\tiny{
\left( \hspace{-0.2cm} \begin{array}{c}
0\\
b_1+b_2 \\
b_3 \\
. \\
.\\
.
\end{array}
 \hspace{-0.2cm}
\right)}}$ s.t.  $\gamma \sim_{X_b} \gamma '$ 
(the assumption of Prop. \ref{step1} 
in our case  is true by Lemma \ref{HOc-ei} and by our assumption in 
($\star$)).  
We have $\tilde{H}_{p-1} (\Delta_{
\tiny{
\left( \hspace{-0.2cm} \begin{array}{c}
0\\
b_1+b_2 \\
b_3 \\
. \\
.\\
.
\end{array}
 \hspace{-0.2cm}
\right)}})= \tilde{H}_{p-1} (\Delta_{
\tiny{
\left( \hspace{-0.2cm} \begin{array}{c}
b_1+b_2 \\
b_3 \\
. \\
.\\
.
\end{array}
 \hspace{-0.2cm}
\right)}})
 =0$, where the last equality holds since 
 ${\cal O}_{{\bf P}^{n_1 - 1}  \times {\bf P}^{n_2} 
 \times .... \times {\bf P}^{n_d}
}(1,..., 1)$  satisfies $N_{p}$ by induction hypothesis.
Thus   $\gamma' \sim 0 $ in $\Delta_{
\tiny{
\left( \hspace{-0.2cm} \begin{array}{c}
0\\
b_1+b_2 \\
b_3 \\
. \\
.\\
.
\end{array}
 \hspace{-0.2cm}
\right)}}$.  
Thus  $\gamma \sim_{X_b} 0 $ and then
 $\gamma \sim_{\Delta_b} 0$
by Prop. \ref{step2}.

Since ${\cal O}_{{\bf P}^{n_1} \times .... \times {\bf P}^{n_d}
}(1,..., 1)$  satisfies $N_{1}$ then, by $(\star) $,
 it satisfies $N_2$ and by applying 
again $(\star)$ we get that it satisfies also $N_3$.

Finally  by Theorem \ref{LPW} and Proposition 
\ref{prodotti} we know that ${\cal O}_{{\bf P}^{n_1} \times .... \times 
{\bf P}^{n_d}
}(1,..., 1)$  does not  satisfy $N_{4}$ if there is $n_i \geq 2$. If $n_i =1$ 
$\forall i$ then we already know the result from Theorem \ref{P1}.
\hfill \framebox(7,7)

\bigskip

\bigskip

Now we will prove  Propositions
 \ref{step1} and  \ref{step2} and  Lemma \ref{HOc-ei}.

\begin{notation}
Let $b \in {\bf N} A_{n_1,..., n_d}$
Let $\gamma$ be a $(p-1)$-cycle in $X_b$.

 For every vertex $a$  in $\gamma$, let
${\cal S}_{a,\gamma} $ be  the set of simplexes of $\gamma$ with vertex  $a$ 
and $\mu_{a, \gamma}$ be the $(p-2)$-cycle s.t. $a \ast \mu_{a, \gamma} 
= \sum_{\tau \in {\cal S}_{a,\gamma}} \tau$.  
For $\tilde{a} \in A_{n_1,..., n_d}$, let 
$$\alpha_{a, {\tilde a}, \gamma} = (a- \tilde{a}) \ast  \mu_{a, \gamma}$$
\end{notation}

\bigskip
 
{\it Proof of Prop. \ref{step1}}.
We order in some way the (finite) vertices of  
$\gamma$ with first coordinate $ \neq 0$: $a^1,..., a^r$.
Let $\tilde{a}^{j}= \tiny{
\left( \hspace{-0.2cm} \begin{array}{c}
0 \\
a^{j}_2 +a^{j}_1\\
a^{j}_3 \\
a^{j}_4 \\
.\\
.
\end{array}
 \hspace{-0.2cm}
\right)}$ for $j=1,..., r$.

Obviously  $\alpha_{a^1,\tilde{a}^1 , \gamma} \sim_{X_b} 0 $,
because  $\mu_{a^1, \gamma} $
 is in $\Delta_{b-a^1}$ and 
$\tilde{H}_{p-1}((a^1 - \tilde{a}^1) \ast  \Delta_{b-a^1}) =
\tilde{H}_{p-2}(\Delta_{b-a^1}) =0$ (since    
${\cal O}_{{\bf P}^{n_1} \times .... \times {\bf P}^{n_d}
}(1,..., 1)$ satisfies Property $N_{p-1}$).
Thus $\gamma_1 := \gamma + \alpha_{a^1, \tilde{a^1},
 \gamma}$ is homologous to $\gamma$ in $X_b$.

We define by induction $\gamma_j := \gamma_{j-1} + \alpha_{a^j,
\tilde{a}^j,  \gamma_{j-1}}$ for $j=2,..., r$.
We want to  prove $\gamma_r \sim_{X_b} 0$; to  prove this, we
 prove  $\alpha_{a^j, \tilde{a^j}, \gamma_{j-1}} \sim_{X_b} 0$
for $j=2,..., r$.

Observe that 
 $\mu_{a^j, \gamma_{j-1} } $ is in
$$\Delta_{b-a^j} \cup \Delta_{
\tiny{
\left( \hspace{-0.2cm} \begin{array}{c}
(b-a^j)_1 -1\\
(b-a^j)_2 +1 \\
(b-a^j)_3 \\
. \\
.\\
.
\end{array}
 \hspace{-0.2cm}
\right)}} \cup .... \cup \Delta_{
\tiny{
\left( \hspace{-0.2cm} \begin{array}{c}
0\\
(b-a^j)_1+(b-a^j)_2 \\
(b-a^j)_3 \\
. \\
.\\
.
\end{array}
 \hspace{-0.2cm}
\right)}}$$

We can find some cycles $\theta_{\varepsilon}$ 
in  $\Delta_{b-a^j - \varepsilon e_1 + \varepsilon e_2}$ for 
 $\varepsilon \in \{0, ..., (b-a^j)_1\}$ s.t. 
 $\mu_{a^j, \gamma_{j-1} }= \sum_{\varepsilon \in \{0, ..., (b-a^j)_1\}}
 \theta_{\varepsilon}$, in fact: if $p=2$ the statement is obvious; let $p 
\geq 3$;  
let $\sigma_0$ be the sum of the simplexes of $\mu_{a^j, \gamma_{j-1} }$ in 
$\Delta_{b-a^j}$ and not in   $\Delta_{b-a^j - e_1 }$; 
$\partial \sigma_0 $ is in  $\Delta_{b-a^j -  e_1 }$ and 
since $\tilde{H}_{p-3}(\Delta_{c- e_1}) =0 $
  $\forall c $ with $ deg(c) \geq p+1$ then
$\exists \sigma_0' $ in $\Delta_{b-a^j -  e_1 }$ s.t. 
$\partial \sigma_0' = \partial \sigma_0 $; 
let $ \theta_0= \sigma_0 - \sigma_0'$; now 
$\mu_{a^j, \gamma_{j-1} } - \theta_0$ is in $$ \Delta_{\tiny{
\left( \hspace{-0.2cm} \begin{array}{c}
(b-a^j)_1 -1\\
(b-a^j)_2 +1 \\
(b-a^j)_3 \\
. \\
.\\
.
\end{array}
 \hspace{-0.2cm}
\right)}} \cup .... \cup \Delta_{
\tiny{
\left( \hspace{-0.2cm} \begin{array}{c}
0\\
(b-a^j)_1+(b-a^j)_2 \\
(b-a^j)_3 \\
. \\
.\\
.
\end{array}
 \hspace{-0.2cm}
\right)}}$$ 
and we can go on analogously: let  $\sigma_1$ 
be the sum of the simplexes of $\mu_{a^j, \gamma_{j-1} }- \theta_0$ in 
$\Delta_{b-a^j-e_1 +e_2}$ and not in   $\Delta_{b-a^j - 2 e_1 }$....

Since $deg(b-a^j -\varepsilon e_1 + \varepsilon e_2 ) \geq p+1$
 and   ${\cal O}_{{\bf P}^{n_1} \times .... \times {\bf P}^{n_d}
}(1,..., 1)$ satisfies $N_{p-1}$ we have
$\tilde{H}_{p-2} (\Delta_{b-a^j -\varepsilon e_1 + \varepsilon e_2} )=0$   thus
$ \theta_{\varepsilon} \sim 0$ in  $\Delta_{
b-a^j -\varepsilon e_1 + \varepsilon e_2}$ 
and then $(a^j -\tilde{a^j}) \ast \theta_{\varepsilon} \sim_{X_b} 0$;
therefore $\alpha_{a^j, \tilde{a}^j, \gamma_{j-1}} \sim_{X_b} 0$.

Thus we can define $\gamma' = \gamma_r$: 
 $\gamma' $ is  in  $\Delta_{
\tiny{
\left( \hspace{-0.2cm} \begin{array}{c}
0\\
b_1+b_2 \\
b_3 \\
. \\
.\\
.
\end{array}
 \hspace{-0.2cm}
\right)}}$ and $ \gamma' \sim_{X_b} \gamma$.
\hfill \framebox(7,7)

\bigskip

{\it Proof of Lemma \ref{HOc-ei}.}  Let $P, Q \in \Delta_{c-e_i}$. Then we can
 find $s \in {\bf N} A_{n_1,..., n_d}$ with $deg (s) \geq 3 $ s.t. $s_i 
\leq c_i $  $\forall i $ and $P, Q \in \Delta_s$; since ${\cal O}_{
{\bf P}^{n_1}  \times.... \times {\bf P}^{n_d}  } (1,..., 1) $ satisfies 
Property $N_1 $   we have that $\tilde{H}_0 (\Delta_s) $; thus there is a 
$1$-cycle $\eta$ in $\Delta_s $ and then in $\Delta_{c-e_i} $ s.t. 
$\partial \eta =P-Q$.
\hfill \framebox(7,7)

\bigskip
Now we will prove Prop. \ref{step2}.

\begin{definition} \label{UFO} Let $d,n_1,...,n_d ,k \in {\bf N} -\{0\}$ and
 $\beta \in {\bf N} A_{n_1,..., n_d}$.
We say that a  $(k-1)$-chain   $\eta$  in $\Delta_{\beta}$ 
  is a $UFO$ with axis  $\langle a^1,..., a^t \rangle $ for the
 coordinate $i$ (for short we will write   
$UFO_{t,k}^{i}(a^1,..., a^t, \Delta_{\beta})$)
if  
$$\eta = \langle a^1,..., a^t \rangle  \ast C_{\eta}$$ for some 
$ (k-t -1)$-cycle $C_{\eta}$ and
 $a^1,..., a^t $ are distinct vertices in $  \Delta_{\beta}$
 with $$(a^1 + ...+ a^t)_i = \beta_i  \;\;\;\;\;\;\;
(a^j)_i > 0 \; \; \forall j =1,...,l$$
We will denote the axis 
$\langle a^1 ,..., a^t \rangle$ by $\chi_{\eta}$.  Observe 
$\partial \eta \subset \Delta_{\beta -e_i}$.
 
(Sometimes we will omit some index when it will be obvious.)
\end{definition}

\smallskip

\hspace{4.5cm}
\includegraphics[scale=0.3]{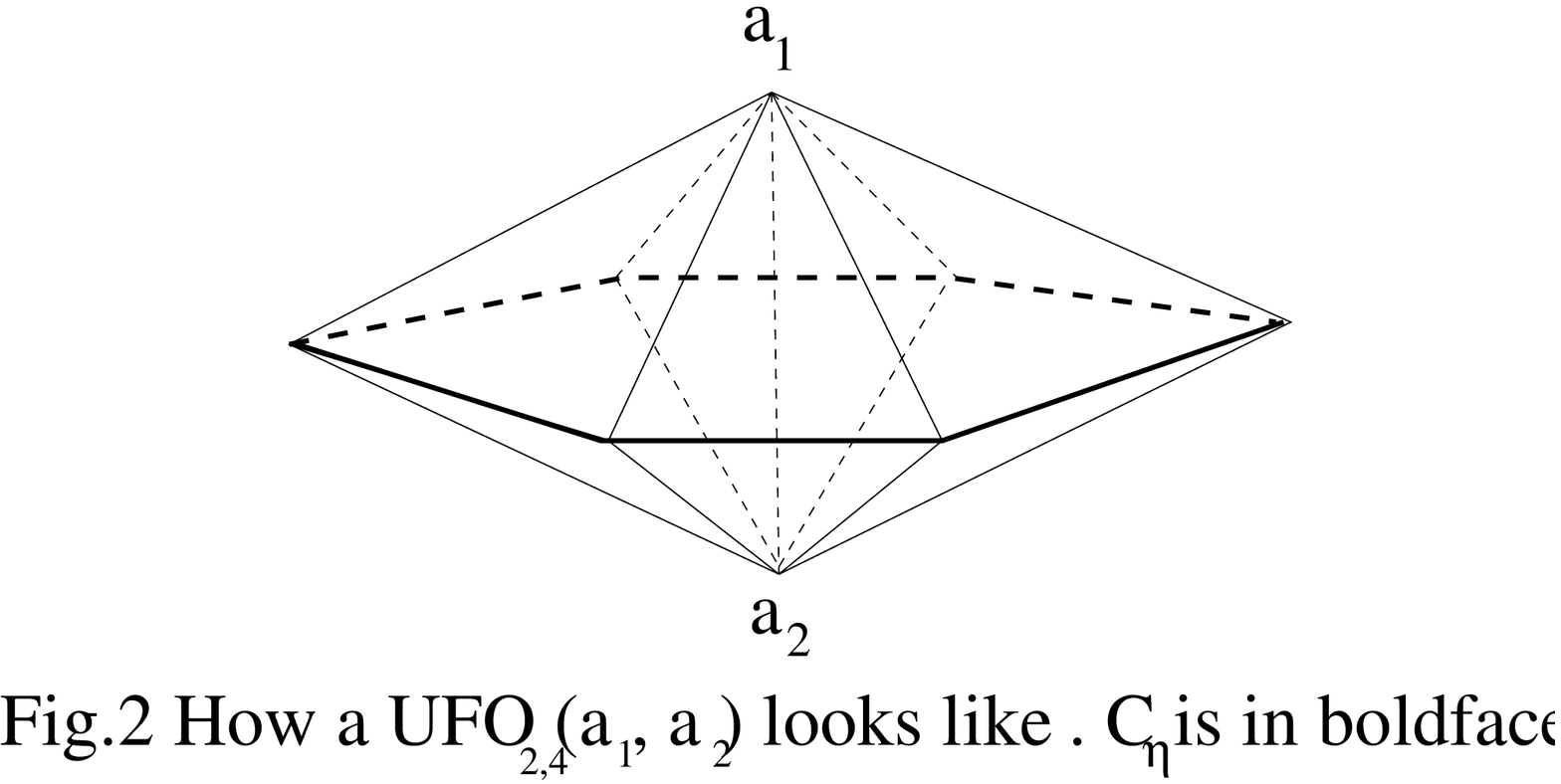}


\begin{lemma} \label{ufisemplici} {\boldmath
 ``$UFO_{p+1, p+1}$, `$UFO_{p, p+1}$, 
`$UFO_{1, p+1}$ ''}
 Let $d, n_1,...n_d, p \in {\bf N}- \{0\}$
 and  $ r, l \in \{1,..., \sum_{i =1,..,d} (n_i +1)\}$ with   $r \neq l$.
 Let $\beta \in {\bf N} A_{n_1,..., n_d}$ with  $\deg \beta \geq p+2$.
Let $\eta $ be  a $UFO_{t, p+1}^r (a^1, ..., a^{t},\Delta_{\beta})$. 

 If $t \in \{p+1, p , 1\}$ 
then $\exists \tilde{\eta} $ $p$-chain  in $\Delta_{
\beta +e_l -e_r}$ with $\partial \tilde{\eta}= \partial \eta$.

\end{lemma}

{\it Proof.} \underline{Case $t=p+1$}
 Since $\deg \beta \geq p+2 $ and $\eta $ is a simplex with 
$p+1$ vertices (the simplex $\langle a^1,..., a^{p+1} \rangle$), then     
 $\exists x \in A_{n_1,..., n_d}$  s.t. $ x \ast \eta  \subset \Delta_{\beta}$.
 Since $(a^1+ ...+a^{p+1})_r = \beta_r$ then  $x_r=0$. Take  $\tilde{\eta} :=
x \ast \partial \eta $.

\underline{Case $t=p$}
 If $\deg \beta \geq p+3$ then $ \deg (\beta -a^1 ...- a^p ) \geq 3$, 
therefore $ \tilde{H}_0 (\Delta_{\beta -a^1 ...- a^p}) =0$, thus  
$\exists \gamma $ s.t. $ \partial \gamma = C_{\eta}$
 and we can take $ \tilde{\eta}=
\gamma \ast \partial  \chi_{\eta}$.

Thus we can suppose $\deg \beta = p+2$. 
Since $C_{\eta}$ is a $0$-cycle it is sufficient to prove the statement when
$C_{\eta}=P-Q$ for some $P,Q \in A_{n_1,..., n_d}$ with
$P=Q +e_i -e_j$  for some $i$ and $j$ (in fact  we can write $C_{\eta}$ as 
$ \sum_s (-1)^s P_s$ with $P_{s+1} $ obtained from $P_s$  by adding 
$1$ to a  coordinate   and subtracting $1$ to another coordinate). 
Let $x = \beta -a^1...-a^p -P$ and  $y =\beta -a^1...-a^p -Q$.
Since $(a^1 + ...+ a^p)_r = \beta_r$ we  have  $x_r =y_r =0$.

Suppose first  $i, j \not \in \{l,r\} $.
Let $z = x +e_l -e_j$.
The chain  $\tilde{\eta } = z  \ast \partial \eta$ is in $ \Delta_{\beta +e_l
-e_r}$ and $\partial \tilde{\eta} = \partial \eta$.

Suppose now   $j \not \in \{l ,r\}$  and $i=l$.
Then  $\tilde{\eta } := x  \ast \partial \eta$ is in $ \Delta_{
\beta +e_l -e_r}$ and $\partial \tilde{\eta} = \partial \eta$.

\underline{Case $t=1$} Let $ \tilde{a^1} = a^1 +e_l -e_r$;
take  $\tilde{\eta} = \tilde{a^1}  \ast C_{\eta }$.
 \hfill \framebox(7,7)

\begin{lemma} \label{subc}  Let $d,n_1, ...,n_d, p \in {\bf N}- \{0\}$ and 
$ r  \in \{1,..., \sum_{i =1,..,d} (n_i +1)\}$ 
Let $\beta \in {\bf N} A_{n_1,..., n_d}$ and $\deg \beta \geq  p+2$.  
Let $\eta $ be  a $UFO_{t, p+1}^r (a^1,..., a^t, \Delta_{\beta })$. 
If $C_{\eta} = \partial \sigma$ where $\sigma $ is a simplex 
(with $p+2 -t $   vertices) in
 $\Delta_{\beta - a^1 -... -a^t} $, then $\exists \tilde{\eta}$ $p$-chain 
in $\Delta_{\beta  -e_r}$ with $\partial \tilde{\eta}= \partial \eta$.
\end{lemma}

{\it Proof.} Since $\eta = \chi_{\eta}   \ast
C_{\eta} =  \chi_{\eta}  \ast   \partial \sigma $, then
 $\partial \eta =  \partial \chi_{\eta}  * \partial \sigma = \partial (
\partial \chi_{\eta}  *  \sigma)$. Take $\tilde{\eta} 
= \partial \chi_{\eta}  *  \sigma$. 

\hfill \framebox(7,7)

\smallskip

\begin{lemma} \label{ufo2}  {\boldmath ``$UFO_{2,4}$''}  
 Let $d,n_1, ...,n_d, p \in {\bf N}- \{0\}$.
Let $\beta \in {\bf N}  A_{n_1,..., n_d}$ with  $\deg \beta \geq 5$.
Let $\eta $ be a $UFO_{2, 4}^2 (a^1, a^2, \Delta_{\beta})$.
 If $n_1 \leq  3 $ then  
 $\exists \tilde{\eta}$ $3$-chain in $\Delta_{\beta
+e_1 -e_2}$ with $\partial \tilde{\eta}= \partial \eta$.
\end{lemma}

{\it Proof.} 
Observe  $C_{\eta}$ is a  $1$-cycle  and obviously $\beta_2 =2$.

$\bullet $ {\em  Case where
 at least one of the first $n_1 +1  $ coordinates of 
$\beta$  except the second one is $\geq 2$,    say $\beta_i  \geq 2$.}

By Lemma \ref{subc} we can substitute a simplex of $C_{\eta} $, $\langle 
T_1,  T_2 \rangle$ with $\langle T_1, P\rangle + \langle P, T_1 \rangle$ where
$P $ is s.t. $\langle P, T_1 , T_2 \rangle \ast \langle a^1 ,a^2 \rangle 
\subset \Delta_{\beta} $ and then we can suppose 
 $C_{\eta }$ is in $ \Delta_{\beta -a^1 -a^2 -e_i}$ 
(if $(T_1+ T_2)_i= (\beta -a^1 -a^2)_i$ we can replace  $\langle 
T_1,  T_2 \rangle$ with $\langle T_1 , P\rangle + \langle P , 
T_1 \rangle$ where
$P $ is s.t. $\langle P, T_1,  T_2 \rangle \ast \langle a^1 ,a^2 \rangle 
\subset \Delta_{\beta} $).
 Thus $\eta $ is  $\Delta_{\beta -e_i}$.

Let  $\lambda := \langle a^2 ,\tilde{a^1} \rangle  \ast C_{\eta}$, 
where $\tilde{a^1}= a^1 - e_2  + e_1 $; $\lambda $ is in $\Delta_{\beta -e_2
 +e_1}$; thus $ \partial \eta  $ is homologous to $\partial (\eta +\lambda)
$ in $  \Delta_{\beta -e_2
 +e_1}$; but $\partial (\eta +\lambda) = (a^1 - \tilde{a^1})  \ast C_{\eta} =
 (a^1 - \tilde{a^2} + \tilde{a^2} - \tilde{a^1})  \ast C_{\eta} = \partial (
(\langle a^1 ,\tilde{a^2}  \rangle + \langle \tilde{a^2} ,
 \tilde{a^1} \rangle )  \ast C_{\eta})$ 
 where 
$\tilde{a^2} = a^2 - e_2 + e_i$; since $
(\langle a^1 ,\tilde{a^2}  \rangle + \langle \tilde{a^2} ,
 \tilde{a^1} \rangle)   \ast C_{\eta}
\subset \Delta_{\beta -e_2+e_1}$, we conclude.

$\bullet $ {\em Case where the first $n_1 +1  $ coordinates of $\beta$  except 
the second one are $\leq 1$.}

In this case, since $n_1 \leq 3$ and $deg (\beta) \geq 5$, we have that 
all the first $n_1 +1  $ coordinates except 
the second one are equal to $1$ (and $n_1=3 $ and $deg(\beta)=5$).

 Let  $\lambda := \langle a^2 ,\tilde{a^1}  \rangle  \ast C_{\eta}$, 
where $\tilde{a^1}= a^1 - e_2  + e_1 $; $\lambda $ is in $\Delta_{\beta -e_2
 +e_1}$; thus $ \partial \eta  $ is homologous to $\partial (\eta +\lambda)
$ in $  \Delta_{\beta -e_2
 +e_1}$; but $\partial (\eta +\lambda) = (a^1 - \tilde{a^1})  \ast C_{\eta} =
 (a^1 - \tilde{a^2} + \tilde{a^2} - \tilde{a^1})  \ast C_{\eta} = \partial (
(\langle a^1 ,\tilde{a^2}  \rangle + \langle \tilde{a^2} ,
 \tilde{a^1} \rangle)   \ast C_{\eta})$ 
 where 
$\tilde{a^2} = a^2 - e_2 + e_1$. 

If every vertex of $C_{\eta}$ has first coordinate $0$, then $ \eta \subset 
\Delta_{\beta -e_1} $ and $
(\langle a^1 ,\tilde{a^2}  \rangle + \langle \tilde{a^2} ,
 \tilde{a^1} \rangle)   \ast C_{\eta}
\subset \Delta_{\beta -e_2+e_1}$ and then we conclude.  

Suppose  there is a  vertex $V$ of $C_{\eta}$ with
 first coordinate equal to $1$ (observe that then, 
 if $\langle P,V \rangle$ is a simplex in $C_{\eta}$, then $P_1 =0$);
 then the simplexes of 
$\langle a^1 ,\tilde{a^2}  \rangle + \langle \tilde{a^2} ,
 \tilde{a^1} \rangle)   \ast C_{\eta}$ not in $ \Delta_{\beta -e_2+e_1}$
(which are the simplexes of 
$\langle \tilde{a^2} ,
 \tilde{a^1} \rangle   \ast C_{\eta}$ not in $ \Delta_{\beta -e_2+e_1}$)
form $UFO^1_{3,4} ( \Delta_{\beta -2e_2+2e_1})$; thus 
by Lemma \ref{ufisemplici} Case $t=p$
we can prove that  $
\partial ((\langle a^1 ,\tilde{a^2}  \rangle + \langle \tilde{a^2} ,
 \tilde{a^1} \rangle)   \ast C_{\eta}) \sim 0$
in $ \Delta_{\beta -e_2+e_1}$ and this finishes the proof.

\hfill \framebox(7,7)

\begin{corollary} \label{fine}  Let $d,n_1, ...,n_d, p \in {\bf N}- \{0\}$ 
with $n_1 \leq 3$
and $p \in \{2,3\}$. Let $\beta \in  {\bf N} A_{n_1,..., n_d}$  and 
$\deg \beta \geq p+2$.
If $\eta$ 
is a $p$-chain  in $\Delta_{\beta}$ with 
 $\partial \eta $ in $\Delta_{
\beta  -e_2}$,
then  $\exists \tilde{\eta}$ $p$-chain in $\Delta_{
\beta +e_1 -e_2}$ with $\partial \tilde{\eta}= \partial \eta$.
\end{corollary}

{\it Proof.}
To prove the statement,  is sufficient to prove it when 
 $\eta$ is a $UFO^2_{t,p+1}$ for $t=1,...,p+1 $, $p \in \{2,3\}$
since $\eta$ is a sum  of $UFO^2_{t,p+1}$ for $t=1,...,p+1 $ $p \in \{2,3\}$.
    Thus our statement  follows from Lemmas \ref{ufisemplici} and \ref{ufo2}.
\hfill \framebox(7,7)

\bigskip

{\it Proof of Prop. \ref{step2}.}
We will show that if $\gamma$ is a $(p-1)$-cycle  in $\Delta_{b} $ with 
$\gamma = \partial \eta$ with  $\eta $ 
 in $\Delta_{b}  \cup \Delta_{b -e_1 + e_2} \cup ...
 \cup \Delta_{b -ke_1 +k e_2}$  for some $k \leq  b_1$, 
 then 
we can construct $\eta'$ in  $\Delta_{b}  \cup ... \cup 
\Delta_{b -(k-1) e_1 +(k-1) e_2}$ s.t. $\partial \eta' = \gamma$
(this, by induction on $b_1$, implies obviously Prop. \ref{step2}):

 let $\nu$ be the 
sum of the simplexes of $\eta$ in 
$\Delta_{b -ke_1 +k e_2 }$
 and not in $\Delta_{b -k e_1 +(k-1) e_2}$;
 $\partial \nu $ is in $\Delta_{
b -k e_1 +(k-1) e_2}$; by Corollary 
\ref{fine} $\partial \nu  = \partial \nu' $ for some 
$\nu'$ in $\Delta_{b -(k-1) e_1 +(k-1) e_2}$
 let  $\eta'= \eta- \nu + \nu' $; $\eta'$ is in 
$\Delta_{b}  \cup ... \cup 
\Delta_{b -(k-1) e_1 +(k-1) e_2}$
 and $\partial \eta' = \partial \eta= \gamma$.
\hfill \framebox(7,7)

\smallskip

{\bf Address: Elena Rubei, Dipartimento di Matematica ``U. Dini'',
via Morgagni 67/A, 
50134 Firenze, Italia.}
{\bf E-mail address: rubei@math.unifi.it}

\smallskip

{\bf 2000 Mathematical Subject Classification:} 
14M25, 13D02.

{\bf Key words:}  syzyzgies, Segre embeddings.

\end{document}